\documentclass[final]{siamltex}

\usepackage{amsmath,amssymb,amsfonts,epsfig}

\title{Linear chord diagrams on two intervals}

\author{J.\ E.\ Andersen$^{1}$,  R.\ C.\  Penner$^{1,2}$, C.\ M.\ Reidys$^{3}$, R.\ R.\ Wang$^{3}$\\[6pt]
$^{1}$ Center for the Quantum Geometry of Moduli Spaces\\
 Aarhus University, DK-8000 Aarhus C, Denmark\\[5pt]
$^{2}$ Departments of Math and Physics\\
 Caltech, Pasadena, CA 91125 USA\\[5pt]
$^{3}$ Center for Combinatorics\\
 LPMC-TJKLC, Nankai University, Tianjin 300071, P.R.~China\\[5pt]}

\begin{document}
\maketitle

{\renewcommand{\thefootnote}{}
\footnote{\emph{E-mail addresses}: andersen{\char'100}imf.au.dk (J.\ E.\ Andersen), rpenner{\char'100}imf.au.dk (R.\ C.\  Penner), duck{\char'100}santafe.edu(C.\ M.\ Reidys), wangrui{\char'100}cfc.nankai.edu.cn (R. R. Wang)}

\footnotetext[1]{JEA and RCP
supported by QGM (Center for Quantum Geometry of Moduli
Spaces, funded by the Danish National Research Foundation).}

\footnotetext[2]{CMR and RRW supported by the 973 Project, the PCSIRT Project of the Ministry of Education, the Ministry of Science and Technology, and the National Science Foundation of China.}

\begin{abstract}
Consider all possible
ways of attaching disjoint chords to two ordered and oriented disjoint intervals so as to produce a connected
graph.  Taking the intervals to lie in the real axis with the induced orientation and the chords to lie in the
upper half plane canonically determines a corresponding fatgraph which
has some associated genus $g\geq 0$, and we consider the natural
generating function ${\bf C}_g^{[2]}(z)=\sum_{n\geq 0} {\bf c}^{[2]}_g(n)z^n$
for the number ${\bf c}^{[2]}_g(n)$ of distinct such chord diagrams of
fixed genus $g\geq 0$ with a given number $n\geq 0$ of chords.
We prove here the surprising fact that ${\bf C}^{[2]}_g(z)=z^{2g+1}  R_g^{[2]}(z)/(1-4z)^{3g+2} $
is a rational function, for $g\geq 0$, where the
polynomial $R^{[2]}_g(z)$ with degree at most $g$ has integer coefficients and satisfies $R_g^{[2]}({1\over 4})\neq 0$.
Earlier work had already determined that the analogous generating function
${\bf C}_g(z)=z^{2g}R_g(z)/(1-4z)^{3g-{1\over 2}}$
for chords attached to a single interval is
algebraic, for $g\geq 1$, where the polynomial $R_g(z)$  with degree at most $g-1$ has integer coefficients and satisfies $R_g(1/4)\neq 0$ in analogy to the generating function ${\bf C}_0(z)$ for the Catalan numbers.  The new results here on ${\bf C}_g^{[2]}(z)$ rely
on this earlier work, and indeed, we find
that $R_g^{[2]}(z)=R_{g+1}(z) -z\sum_{g_1=1}^g R_{g_1}(z) R_{g+1-g_1}(z)$, for $g\geq 1$.
\end{abstract}



\pagestyle{myheadings}
\thispagestyle{plain}
\markboth{J. E. ANDERSEN, R. C. PENNER, C. M. REIDYS, AND R. R. WANG}{Linear chord diagrams on two intervals}


\section{Introduction}\label{S:intro}


Fix a collection of disjoint oriented intervals called ``backbones'' in a specified
linear ordering and consider all possible ways of attaching a family of
unoriented and unordered intervals called ``chords'' to the backbones by gluing
their endpoints to pairwise distinct points in the backbone.  These combinatorial structures
occur in a number of instances in pure mathematics including finite type invariants
of links \cite{Barnatan95,Kontsevich93}, the geometry of moduli spaces of flat connections
on surfaces \cite{AMR1,AMR2}, the representation theory of Lie algebras \cite{CSM}, and the mapping class groups  \cite{BAMP}.  They also
arise in applied mathematics for codifying the
pairings among nucleotides in a collection of RNA molecules \cite{Reidys},
or more generally in any macromolecular complex \cite{Nagai-Mattaj},
and
in bioinformatics where they apparently
represent the building blocks of grammars of folding algorithms
of RNA complexes \cite{Alkan,backofen,RIP}.

This paper is dedicated to enumerative problems associated with
connected chord diagrams on two backbones as follows.

Drawing a chord diagram $G$ in the plane with its backbone segments
lying in the real axis with the natural orientation and its chords in the upper half-plane determines
cyclic orderings on the half-edges of the underlying graph
incident on each vertex.  This defines a corresponding ``fatgraph''
$\mathbb G$, to which is canonically associated a topological surface
$F(\mathbb G)$ (cf.\ $\S$\ref{S:fatgraphs}) of some genus, cf.\ Figure~\ref{F:two}.

\begin{figure}[ht]
\centerline{\epsfig{file=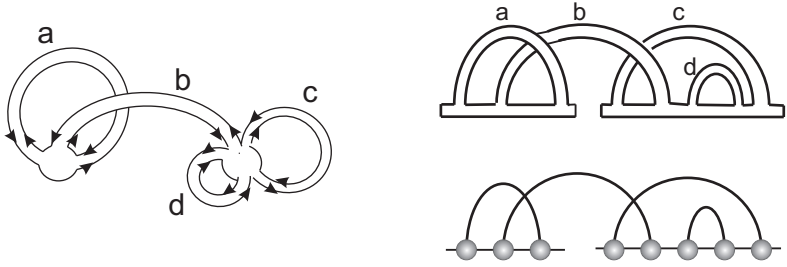,width=0.8\textwidth} \hskip8pt}
\caption{\small From a chord diagram to a fatgraph and to a surface.}
\label{F:two}
\end{figure}

Let ${\bf c}_g^{[b]}(n)$ denote the number of distinct
chord diagrams on $b\geq 1$ backbones with $n\geq 0$ chords
of genus $g\geq 0$ with its natural generating function
${\bf C}_g^{[b]}(z)=\sum_{n\geq 0} {\bf c}^{[b]}_g(n)z^n$,
setting in particular ${\bf c}_g(n)={\bf c}_g^{[1]}(n)$ and
likewise ${\bf C}_g(z)={\bf C}_g^{[1]}(z)$ for convenience.
In particular, the Catalan numbers ${\bf c}_0(n)$, i.e., the number of
triangulations of a polygon with $n+2$ sides, enumerate linear chord
diagrams of genus zero, and we have ${\bf C}_0(z)={{1-\sqrt{1-4z}}\over{2z}}$
with

The one-backbone generating functions
\begin{eqnarray*}
\mathbf{C}_g(z) = \, \frac{P_g(z)}{(1-4z)^{3g-{1\over 2}}}, ~{\rm for~any}~g\geq 1,
\end{eqnarray*}
were computed in \cite{APRW},
where $P_g(z)$ is a polynomial defined over the integers of degree at most
$(3g-1)$ that is divisible by $z^{2g}$ with $P_g(1/4)\neq 0$.
In particular, $\mathbf{C}_g(z)$ is algebraic over $\mathbb C(z)$ for all
$g\geq 1$,  just as is the Catalan generating function $\mathbf{C}_0(z)$,
and there are explicit expressions such as
\begin{eqnarray*}
\mathbf{c}_1(n) & = & \frac{2^{n-2}(2n-1)!!}{3(n-2)!}, \\
\mathbf{c}_2(n) & = & \frac{2^{n-4}(5n-2)(2n-1)!!}{90 (n-4)!}, \\
\mathbf{c}_3(n) & = &\frac{2^{n-6} (35n^2-77n+12) (2n-1)!!
                      }{5670 (n-6)! },
\end{eqnarray*}
for these ``higher'' Catalan numbers,
where $\mathbf{c}_g(n)=0$, for $n<2g$.

In Theorem \ref{T:main} of this paper, we derive the two-backbone generating functions
\begin{eqnarray*}
\mathbf{C}_g^{[2]}(z) & = & \frac{P_g^{[2]}(z)}{(1-4z)^{3g+2}},~{\rm for~any}~g\geq 0,
\end{eqnarray*}
where $P_g^{[2]}(z)$ is an integral polynomial of degree at most
$(3g+1)$ that is divisible by $z^{2g+1}$ with $P_g^{[2]}(1/4)> 0$.
In particular, these generating functions are rational functions defined over the integers.
In fact, these polynomials
$$P_g^{[2]}(z) =  z^{-1}P_{g+1}(z)-\sum_{g_1=1}^{g}P_{g_1}(z)P_{g+1-g_1}(z)$$
are computable in terms of the previous ones, for example:
\begin{eqnarray*}
P_0^{[2]}(z)&=&z,\\
P_1^{[2]}(z)&=&z^3(20z+21),\\
P_2^{[2]}(z)&=&z^5\, \left(1696z^2+6096z+1485 \right),\\
P_3^{[2]}(z)&=&z^7\, \left(330560z^3+2614896z^2+1954116z+225225 \right),\\
P_4^{[2]}(z)&=&z^9\left(118652416z^4 +1661701632z^3+2532145536z^2+851296320z+
59520825\right),\\
P_5^{[2]}(z)&=&z^{11} \left(68602726400z^5+
1495077259776z^4\right.\\
&&\left.+3850801696512z^3+2561320295136z^2 +505213089300z+24325703325\right).
\end{eqnarray*}
The experimental fact that the coefficients of $P_g^{[2]}$ are positive, just as are those of
$P_g(z)$, leads us to speculate that they themselves solve an as yet unknown
enumerative problem.
Furthermore, explicit expressions for the number of two-backbone
chord diagrams of fixed genus such as
\begin{eqnarray*}
\mathbf{c}^{[2]}_0(n) & =&  n\, 4^{n-1},\\
\mathbf{c}^{[2]}_1(n) & =&  \frac{1}{12}\left( 13\,n+3 \right)
                          n\left( n-1 \right)\left( n-2 \right)\,{4}^{n-3},\\
\mathbf{c}^{[2]}_2(n) & =&  {\frac {1}{180}}\,\left(
                             445\,{n}^{2}-401\,n-210 \right) n
                       \left(n-1 \right)\left(n-2 \right)\left( n-3 \right)
                                              \left(n-4 \right)\, {4}^{n-6} ,
\end{eqnarray*}
can also be derived, cf.\ Corollary \ref{C:formel}.

\section{Background and Notation}\label{Back}

We formulate the basic terminology for graphs and chord diagrams,
establish notation and recall facts about symmetric groups, recall the
fundamental concepts about fatgraphs, and finally recall and extend
results from \cite{APRW} for application in subsequent sections.

\subsection{Graphs and chord diagrams on several backbones}

A {\it graph} $G$ is a finite one-dimensional CW complex comprised of
one-dimensional open edges $E(G)$ and zero-dimensional vertices $V(G)$.
An edge of the first barycentric subdivision of $G$ is called a {\it half-edge}.  A
half-edge of $G$ is {\it incident} on $v\in V(G)$ if $v$ lies in its closure, and the {\it
valence} of $v$ is the number of incident half-edges.

Consider a collection $\beta_1,\ldots ,\beta_b$ of pairwise disjoint closed intervals with integer endpoints in the real axis $\mathbb R$, where $\beta_i$ lies to the left of $\beta_{i+1}$, for $i=1,\ldots , b-1$; thus,
these {\it backbone intervals} $\{ \beta _i\}_1^b$ are oriented and ordered by the orientation
on $\mathbb R$.  The {\it backbone} $B=\sqcup\{ \beta _i\}_1^b$ itself is regarded as a graph
with vertices given by the integer points $\mathbb Z\cap B$ and edges determined by the unit intervals
with integer endpoints that it contains.  A {\it chord diagram} on the backbone $B$ is a graph
$C$ containing $B$ so that $V(C)=V(B)$ where each vertex in $B-\partial B$ has valence
three and each vertex in $\partial B$ has valence two; in particular, $\# V(G)$ is even,
where $\#X$ denotes the cardinality of the set $X$.
Edges in $E(B)\subseteq E(C)$ are called {\it backbone edges}, and edges
in the complement $E(C)-E(B)$ are called {\it chords}.

\subsection{Permutations}

The symmetric group $S_{2n}$ of all permutations on $2n$ objects plays
a central role in our calculations, and we establish here standard
notation and recall fundamental tools. Let $(i_1,i_2,\ldots ,i_k)$ denote
the cyclic permutation $i_1\mapsto
i_2\mapsto\cdots\mapsto i_k\mapsto i_1$ on distinct objects $i_1,\ldots,
i_k$.  We shall compose permutations $\pi,\tau$ from right to left, so that
$\pi\circ\tau(i)=\pi(\tau(i))$.  Two cycles are {\it disjoint} if they have disjoint supports.

Let $[\pi ]$ denote the conjugacy class of $\pi\in S_{2n}$.  Conjugacy
classes in $S_{2n}$ are
identified with classes of partitions of $\{1,\ldots
,2n\}$, where in the standard slightly abusive notation,
$\pi\in[\pi ]=[1^{\pi_1} \cdots {2n}^{\pi_{2n}}]$ signifies that
$\pi$ is comprised of $\pi_k\geq 0$ many $k$ cycles with pairwise disjoint
supports, for $k=1,\ldots ,2n$.
In particular,  $\sum_{k=1}^{2n} k\pi_k=2n$, and the number of elements in the class
$[\pi]$ is given by
$$
|\pi|= |[\pi]|= {{(2n)!}\over{\prod_{k=1}^{2n}\left( k^{\pi_k}\pi_k!\right)}}.
$$

We shall make use of the fact that the irreducible characters
$\chi^Y$ of $S_{2n}$ are labeled by Young tableaux $Y$ and shall
slightly abuse notation writing $\chi^Y([\pi])=\chi^Y(\pi)$ for the value
taken on either a permutation $\pi$ or its conjugacy class $[\pi]$.
In order to evaluate characters, we shall rely heavily on the
Murnaghan-Nakayama rule
\begin{equation}\label{E:MN}
\chi^Y( (i_1,\dots,i_m ) \circ \sigma ) =
\sum_{Y_\mu~{\rm so~that}~ Y- Y_\mu\, \text{\rm is a} \atop
      \text{\rm rim hook of length $m$}} (-1)^{w(Y_\mu)}
                                          \chi^{Y_\mu}(\sigma),
\end{equation}
where $w(Y_\mu)$ is  one less than the number of rows in the rim hook $Y-Y_\mu$.
See a standard text such as \cite{Sagan} for further details.

\subsection{Fatgraphs}\label{S:fatgraphs}

If $G$ is a graph, then a {\it fattening} of $G$ is the specification of a collection of
cyclic orderings, one such ordering on the half-edges incident on each
vertex of $G$.  A graph $G$ together with a fattening is called a {\it fatgraph}
${\mathbb G}$.

The key point is that a fatgraph $\mathbb G$ determines an oriented surface
$F({\mathbb G})$ with boundary as follows.
For each $v\in V(G)$, take an oriented polygon
$P_v$ with $2k$ sides, and choose a tree in $P_v$ with a univalent
vertex in every other side of $P_v$ joined to a single $k$-valent
vertex in the interior.  Identify the half-edges incident on $v$ with the edges
of the tree in the natural way
so that the induced counter-clockwise cyclic ordering on the boundary of
$P_v$ agrees with the fattening of $\mathbb G$ about $v$.
We construct the surface
$F(\mathbb G)$ as the quotient of $\sqcup_{v\in V(G)}
P_v$ by identifying pairs of frontier edges of the polygons with an orientation-reversing homeomorphism if
the corresponding half-edges lie in a common edge of $G$.  The oriented surface $F(\mathbb G)$
is connected if $G$ is and has some associated genus $g(\mathbb G)\geq 0$ and number $r(\mathbb G)\geq 1$ of boundary components.

Furthermore, the trees in the various polygons combine to give
a graph identified with $G$ embedded in $F(\mathbb G)$, and we may thus regard
$G\subseteq F(\mathbb G)$. By construction, $G$ is a deformation retraction of
$F(\mathbb G)$, hence their  Euler characteristics
agree, namely,
$$
\chi(G)=\#V(G)-\#E(G)=\chi (F(\mathbb G))=2-2g(\mathbb G)-r(\mathbb G),
$$
where the last equality holds provided $G$ is connected.

A fatgraph $\mathbb G$ is uniquely determined by a pair of permutations on
the half-edges of its underlying graph $G$ as follows.  Let
$v_k\geq 0$ denote the number of $k$-valent vertices, for each
$k\geq 1,\ldots ,K$, where $K$ is the maximum valence of vertices of $G$,
so that $\sum_{k=1}^K kv_k=2n$, and
identify the set of half-edges
of $G$ with the set $\{ 1,2,\ldots ,2n\}$ once and for all.  The orbits of the cyclic orderings in
the fattening on $G$ naturally give rise to disjoint cycles comprising a permutation $\tau\in S_{2n}$ with $\tau\in [1^{v_1}2^{v_2}\cdots K^{v_K}]$.  The second permutation $\iota\in [2^n]$ is the product
over all edges of $G$ of the two-cycle permuting the two half-edges it contains.

One important point is that the boundary
components of $F(\mathbb G)$ are in bijective correspondence with the
cycles of $\tau\circ\iota$, i.e.,
$r(\mathbb G)$ is the~number~of~disjoint~cycles~comprising $\tau\circ\iota$.
Furthermore, isomorphism classes of fatgraphs with vertex valencies
$(v_k)_{k=1}^K$ are in one-to-one correspondence with conjugacy classes of
pairs $\tau,\iota\in S_{2n}$, where $\tau\in [1^{v_1}2^{v_2}\cdots K^{v_K}]$
and $\iota\in [2^n]$, cf.\ Proposition \ref{P:isos}.  Thus,
fatgraphs are easily stored and manipulated on the computer, and
various enumerative problems can be formulated in terms of Young tableaux.

See \cite{Penner88,PKWA} for more details on fatgraphs and
\cite{APRW, BIZ,Harer-Zagier,Itzykson-Zuber,Milgram-Penner,Penner92}
for examples of fatgraph enumerative problems in terms of character
theory for the symmetric groups.

\subsection{Fatgraphs and chord diagrams}

A convenient graphical method to uniquely determine a fatgraph is
the specification of a
regular planar projection of a graph in 3-space, where the counter-clockwise cyclic
orderings in the plane of projection determine the fattening and the crossings of edges
in the plane of projection can be arbitrarily resolved into
under/over crossings without affecting the resulting isomorphism class.
A band about each edge can be added to a neighborhood of
the vertex set in the plane of projection respecting orientations in
order to give an explicit representation of the associated surface embedded
in 3-space as on the top-right of Figure \ref{F:two}.

In particular, the standard planar representation of a chord diagram $C$ on two
backbones positions the two backbone segments in the real axis respecting their ordering
and orientation and places the chords as semi-circles in the upper half-plane as on the
bottom right in Figure \ref{F:two}.  This determines the {\it canonical fattening} $\mathbb C$ on $C$ as illustrated on
the top right in the same figure.  We may furthermore collapse each backbone segment to a
distinct point and induce a fattening on the quotient in the natural way to produce a fatgraph $\mathbb G_C$
with two vertices and corresponding surface $F(\mathbb G_C)$ as depicted on the left in Figure \ref{F:two}.
Insofar as $C$ is a deformation retraction of $F(\mathbb C)$ and hence of the homeomorphic surface
$F(\mathbb G_C)$,
these three spaces share the same Euler characteristic.

A chord diagram $C$ on two backbones with $n\geq 0$ chords
therefore gives rise to a
fatgraph $\mathbb G_C$ with two ordered vertices of respective valencies $c$ and
$(2n-c)$ and two distinguished half-edges, namely, the
ones coming just after the locations of the collapsed backbones.
Upon labeling the half-edges, $C$ determines
a pair of permutations $\tau_c\in[c,(2n-c)]$ with its two disjoint cycles corresponding
to the fattening of the two vertices of $\mathbb G_C$ and $\iota\in[2^n]$
corresponding to the edges of $\mathbb G_C$.  Provided $C$ and hence
$\mathbb G_C$ is connected, we may define
$$\aligned
r(C)&=r(\mathbb G_C)=
{\rm the~number~of~disjoint~cycles~comprising}~\tau_c\circ\iota,\\
g(C)&=g(\mathbb G_C)= {1\over 2}\bigl(n-r(C)\bigr).\\
\endaligned$$  Conversely, the specification of $\tau_c\in[c,(2n-c)]$ and
$\iota\in[2^n]$ uniquely determines a chord diagram on two backbones with $2n$
edges.    In fact, the isomorphism class of the fatgraph clearly corresponds to the
conjugacy class of the pair $\tau_c,\iota$.

For example, labeling the half-edges contained in chords from left to right on the chord diagram $C$ on the
bottom-right in Figure ~\ref{F:two} produces the fatgraph on the left of the figure
with permutations given by $\tau_3=(1,3,2)(4,5,6,7,8)$ and $\iota = (1,3)(2,5)(4,8)(6,7)$
with $r(C)=4$ and $g(C)=0$.

Summarizing, we have

\begin{proposition}\label{P:isos}
The following sets are in bijective correspondence:
$$\aligned
&\{ {\rm isomorphism~classes~of~chord~diagrams~on~two~backbones~with}~n~{\rm chords}\}\\
\approx &\{{\rm isomorphism~classes~of~fatgraphs~with~two~vertices~and}~n~{\rm edges}:\\
&\hskip 1.5ex{\rm each~vertex~has~a~distinguished~incident~half-edge}\}\\
\approx&\{{\rm conjugacy~classes~of~pairs}~\tau,\iota\in S_{2n}: \iota\in [2^n]~{\rm and}~ \tau\in[c,(2n-c)]\\
&\hskip 1.5ex{\rm for~some}~1\leq c\leq 2n-1, ~{\rm with~an~ordered~pair~of~elements,~one~from~each~of~the~two~cycles~of}~\tau\}
\endaligned$$
\end{proposition}

Our first counting results will rely upon the bijection established here
between chord diagrams and pairs of permutations.

\subsection{The one backbone case}

We collect here results from \cite{APRW} which will be required in the sequel
as well as extend from \cite{APRW} an explicit computation we shall also need.

As a general notational point for any power series $T(z)=\sum a_iz^i$, we
shall write $[z^i]T(z)=a_i$ for the  extraction of the coefficient $a_i$
of $z^i$.

As mentioned in the introduction, the generating function ${\bf C}_g(z)=\sum_{n\geq 0} {\mathbf c}_g(n)z^n$ for $g\geq 1$ satisfies
$${\bf C}_g(z)=P_g(z)(1-4z)^{{1\over 2}-3g},$$
where the polynomial $P_g(z)$ has integer coefficients, is divisible by $z^{2g}$ and has degree at most
$3g-1$.  Indeed,
$\mathbf{C}_g(z)$ is recursively computed using the ODE
\begin{eqnarray}\label{E:ODE}
z(1-4z)\frac{d}{dz}\mathbf{C}_g(z) +(1-2z)\mathbf{C}_g(z) & = &
\Phi_{g-1}(z),
\end{eqnarray}
where
\begin{eqnarray*}
\Phi_{g-1}(z) =  z^2\left(4z^3\frac{d^3}{dz^3}\mathbf{C}_{g-1}(z) +
24z^2 \frac{d^2}{dz^2}\mathbf{C}_{g-1}(z) + 27z\frac{d}{dz}
\mathbf{C}_{g-1}(z)+3\mathbf{C}_{g-1}(z)\right),
\end{eqnarray*}
with initial condition $\mathbf{C}_g(0)=0$.
This ODE in turn follows from the recursion
\begin{eqnarray}\label{E:recursion}
(n+1)\, \mathbf{c}_g(n) & = &  2(2n-1)\,\mathbf{c}_g(n-1)+
                          (2n-1)(n-1)(2n-3)\,\mathbf{c}_{g-1}(n-2),
\end{eqnarray}
where $\mathbf{c}_g(n)=0$ for $2g>n$,
which is derived from a fundamental identity
first proved in \cite{Harer-Zagier}. Namely, the polynomials
\begin{equation}
P(n,N)=\sum_{\{g:2g\leq n\}}\mathbf{c}_g(n) N^{n+1-2g}
\end{equation}
combine into the generating function
\begin{equation}\label{E:polynN}
1+2\sum_{n\geq 0} \frac{P(n,N)}{(2n-1)!!} z^{n+1}=\biggl ( \frac{1+z}{1-z}\biggr )^N.
\end{equation}
In addition to the explicit expressions for  ``higher'' Catalan number,
one also computes
\begin{equation}\label{E:values}
[z^{2g}]P_g(z)={\bf c}_g(2g)=\frac{(4g)!}{4^g(2g+1)!}.
\end{equation}

In fact, it is also shown in \cite{APRW} that $P_g(1/4)$ is non-zero, but here,
we shall require its exact numerical value for subsequent estimates:

\begin{lemma}\label{L:gamma}We have for $g\geq 1$,
\begin{equation}\label{E:q_g4}
P_g(1/4) =\bigl (\frac{9}{4}\bigr )^g~ {\frac {\Gamma  \left( g-1/6 \right)
           \Gamma  \left( g+1/2\right) \Gamma
           \left( g+1/6 \right) }{6{\pi }^{3/2}
           \Gamma \left( g+1 \right) }}.
\end{equation}
\end{lemma}
\begin{proof}
We compute
$P_g(1/4)$ by induction on $g$, and for the basis step
$g=1$,
\begin{equation}
{\bf C}_1(z)= {\frac {{z}^{2}}{\left( 1-4\,z \right)^{3}}}
              \sqrt {1-4\,z}= {\frac {{z}^{2}}{\left( 1-4\,z \right)^{3}}} \bigl (1-2z{\bf C}_0(z)\bigr )
\end{equation}
according to \cite{APRW}, whence $P_1(1/4)=(1/4)^2=1/16$.

The solution to the ODE (\ref{E:ODE})  is given by
\begin{equation}\label{E:intsoln}
{\bf C}_{g+1}(z)=\, \left(\int_0^z
\frac{\Phi_{g}(y)}{(1-4y)^{\frac{3}{2}}} dy+C \right)\,\frac{\sqrt{1-4z}}{z}=\,\left(\int_0^z
\frac{Q_g(y)}{(1-4y)^{3g+4}} dy+C \right)\,\frac{\sqrt{1-4z}}{z},
\end{equation}
where $Q_g(z)$ is shown to be a polynomial of degree at most $(3g+2)$.

Insofar as ${\bf C}_{g}(z)=P_{g}(z)/(1-4z)^{3g-1/2}$,
we find
$$
\aligned
\frac{d{\bf C}_{g}(z)}{dz}&=
\frac{P_{1g}(z)}{(1-4z)^{3g+1/2}},~{\rm where}~
P_{1g}(z)=(1-4z)P_{g}'(z)+(12g-2)P_{g}(z),\\
\quad \frac{d^2 {\bf
C}_{g}(z)}{dz^2}&=
\frac{P_{2g}(z)}{(1-4z)^{3g+3/2}},~{\rm where}~
P_{2g}(z)=(1-4z)P_{1g}'(z)+(12g+2)P_{1g}(z), \\
\quad
\frac{d^3 {\bf C}_{g}(z)}{dz^3}&=
\frac{P_{3g}(z)}{(1-4z)^{3g+5/2}},~{\rm where}~
P_{3g}(z)=(1-4z)P_{2g}'(z)+(12g+6)P_{2g}(z).\\
\endaligned
$$
Thus,
\begin{equation}\label{E:Q}
Q_g(z)=4z^5P_{3g}(z)+24z^4(1-4z)P_{2g}(z)
+27z^3(1-4z)^2P_{1g}(z)+3z^2(1-4z)^3P_{g}(z),
\end{equation}
whence
\begin{equation}
Q_g(1/4)=4^{-4}P_{3g}(1/4)=4^{-4}(12g+6)(12g+2)(12g-2)P_g(1/4)\neq 0.
\end{equation}
Since $Q_g(z)$ has degree at most $(3g+2)$, its partial fraction expansion
reads
\begin{equation}\label{E:pfrac}
\frac{Q_g(z)}{\left( 1-4\,z \right) ^{(3g+4)}} =
               \sum_{j= 2}^{3g+4}\frac{A_j}{\left( 1-4\,z \right) ^{j}},
\end{equation}
where the $A_j\in \mathbb{Q}$ and $A_{3g+4}=Q_g(1/4)$.
According to \cite{APRW}, $P_{g+1}(z)$ is given by
\begin{equation}\label{E:laurent}
P_{g+1}(z)= {{-1}\over 4z}
\left(\sum_{j= 2}^{3g+4}\frac{-A_j}{j-1}\left( 1-4\,z \right)^{3g+4-j}+
\sum_{j= 2}^{3g+4}\frac{A_j}{j-1}\left( 1-4\,z \right)^{3g+3}\right),
\end{equation}
and hence
\begin{eqnarray}\label{E:recgamma}
P_{g+1}(1/4)&=&4^{-4}(12g+6)(12g+2)(12g-2)P_g(1/4)/(3g+3)\nonumber\\
&=&\frac{9(g+1/2)(g+1/6)(g-1/6)}{4(g+1)}P_g(1/4),
\end{eqnarray}
where $P_1(1/4)=1/16$. The lemma follows by checking
that the formula in eq.~(\ref{E:q_g4}) is the unique solution of
eq.~(\ref{E:recgamma}) with $P_1(\frac14) = \frac{1}{16}$.
\end{proof}


\section{Lemmas on characters}


\begin{lemma}\label{L:kro} For any two permutations $\tau,\pi\in S_{2n}$, we have
\begin{equation}
\sum_{\sigma\in S_{2n}}
\delta_{[\sigma],[2^n]}\cdot\delta_{[\tau\sigma],[\pi]}=
\frac{(2n-1)!!}{\prod_jj^{\pi_j}\cdot \pi_j!}
\sum_{Y} \frac{\chi^Y([2^n])\chi^Y(\pi)\chi^Y(\tau)}
{\chi^Y([1^{2n}])}.
\end{equation}
\end{lemma}
\begin{proof}
In order to prove the lemma, we shall apply orthogonality and completeness
of the irreducible characters of $S_{2n}$, that irreducible
representations are indexed by Young diagrams $Y$ containing $2n$ squares,
and the fact that $\chi^Y(\sigma)=\chi^Y(\sigma^{-1})$, for any $\sigma\in S_{2n}$.

Fix some $\pi\in S_{2n}$ with
$\pi\in[1^{\pi_1} \cdots {2n}^{\pi_{2n}}]$.
According to the orthogonality relation
\begin{eqnarray*}
\sum_Y \chi^Y(\sigma_1)\chi^Y(\sigma_2) & = &
\frac{(2n)!}{\vert [\sigma_1]\vert}\cdot
\delta_{[\sigma_1],[\sigma_2]}
\end{eqnarray*}
of the second kind, we have
\begin{eqnarray*}
\sum_{\sigma\in S_{2n}}
\delta_{[\sigma],[2^n]}\cdot\delta_{[\tau\sigma],[\pi]} & = &
\sum_{\sigma\in S_{2n}}\left[\frac{\vert [2^n]\vert}{(2n)!}
\sum_Y \chi^Y(\sigma)\chi^Y([2^n])\right]
\cdot \left[ \frac{\vert [\pi]\vert}{(2n)!}
\sum_{Y'} \chi^{Y'}(\tau\sigma)\chi^{Y'}(\pi)\right] \\
&=&  \frac{(2n-1)!!}{\prod_j j^{\pi_j}\cdot \pi_j!}
\sum_{Y,Y'} \chi^Y([2^n])\chi^{Y'}(\pi)
\left[\frac{1}{(2n)!}\sum_{\sigma\in S_{2n}}\chi^Y(\sigma)\chi^{Y'}
(\tau\sigma)\right].
\end{eqnarray*}
The variant
\begin{eqnarray*}
\frac{1}{(2n)!}\sum_{\sigma\in S_{2n}}\chi^Y(\sigma)\chi^{Y'}(\tau\sigma)=
\frac{\chi^Y(\tau)}{\chi^Y([1^{2n}])}\cdot \delta_{Y,Y'}
\end{eqnarray*}
of the orthogonality relation of the first kind gives
\begin{eqnarray*}
\sum_{\sigma\in S_{2n}}\delta_{[\sigma],[2^n]}\cdot\delta_{[\tau\sigma],[\pi]}
&=& \frac{(2n-1)!!}{\prod_j j^{\pi_j}\cdot \pi_j!}
\sum_{Y,Y'} \chi^Y([2^n])\chi^{Y'}(\pi)\left[
\frac{\chi^{Y'}(\tau)}{\chi^{Y'}([1^{2n}])}\cdot\delta_{Y,Y'}\right]\\
&=& \frac{(2n-1)!!}{\prod_j j^{\pi_j}\cdot \pi_j!}
\sum_{Y} \frac{\chi^Y([2^n])\chi^Y(\pi)\chi^Y(\tau)}{
\chi^Y([1^{2n}])}
\end{eqnarray*}
as required.
\end{proof}

\begin{figure}[ht]
\centerline{\epsfig{file=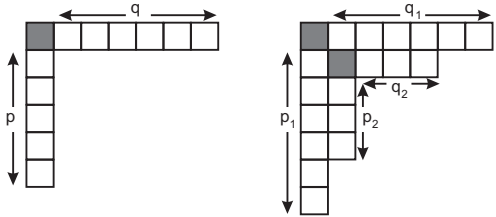,width=0.5\textwidth} \hskip8pt}
\caption{\small The Young tableaux $Y_{p,q}$ and $Y_{p_1,p_2}^{q_1,q_2}$.}
\label{F:Young}
\end{figure}

In the course of our analysis, the two shapes of Young diagrams illustrated in Figure ~\ref{F:Young} are
of importance, where

$\bullet$ $p,q\ge 0$ with $p+q+1=2n$ determine the shape
          $Y_{p,q}$ which is $(p,q)$-hook, having a single row
          of length $q+1\ge 1$ and $p$ rows of length one
          with corresponding character $\chi^{p,q}=\chi^{Y_{p,q}}$,

$\bullet$ $p_1\ge p_2+1\ge 1$, $q_1\ge q_2+1\ge 1$ with
          $p_1+q_1+p_2+q_2+2=2n$ determine the shape
           $Y_{p_1,p_2}^{q_1,q_2}$ with
          one row of length $q_1+1$, one row of length $q_2+2$,
          $p_2$ rows of length two and $p_1-p_2-1$ rows of length
          one.\\

\begin{lemma}\label{L:char}
Suppose $\tau_c=(1,\dots,c)\circ (c+1,\dots,2n)$.
Then we have
\begin{equation}\label{E:1L}
\sum_{c=1}^{2n-1} \chi^{Y}(\tau_c)= (-1)^p(q-p)\,\delta_{Y,Y_{p,q}}
\end{equation}
for any Young diagram $Y$.
\end{lemma}
\begin{proof}
\begin{figure}[ht]
\centerline{\epsfig{file=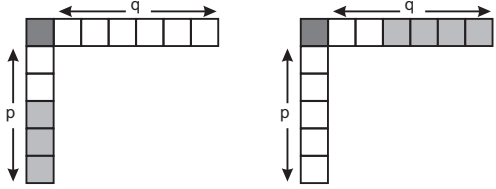,width=0.5\textwidth} \hskip8pt}
\caption{\small Rim-hooks in a Young diagram $Y_{p,q}$ are either vertical
or horizontal blocks. }
\label{F:Youngrim1}
\end{figure}
In order to prove eq.~(\ref{E:1L}), let us first assume $Y=Y_{p,q}$.
Since the only rim-hooks in $Y$ lie entirely in the first row or column as in Figure
\ref{F:Youngrim1}
the Murnaghan-Nakayama rule gives
\begin{eqnarray*}
\sum_{c=1}^{2n-1}\chi^{p,q}(\tau_c)
  & = & \sum_{c=1}^{2n-1}\left(\chi^{p,q-c}((c+1,\dots,2n)) +
                         (-1)^{c-1}\, \chi^{p-c,q}((c+1,\dots,2n)) \right) \\
  & = & \sum_{c=1}^{q}(-1)^p + \sum_{c=1}^{p} (-1)^{c-1}(-1)^{p-c}  \\
                   & = & q(-1)^p + p(-1)^{p-1}\\
                   &=& (-1)^{p}(q-p).
\end{eqnarray*}

\begin{figure}[ht]
\centerline{\epsfig{file=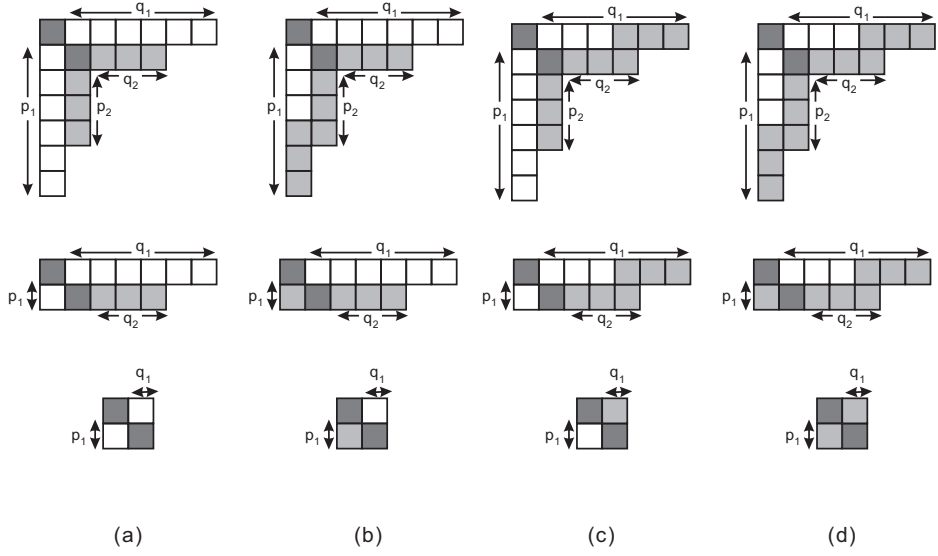,width=0.95\textwidth}
\hskip8pt} \caption{\small The case of $Y_{p_1,p_2}^{q_1,q_2}$: for
the character value $\chi^{Y_{p_1,p_2}^{q_1,q_2}}$ a rim-hook of
size $c$ and a hook of size $2n-c$ has to be removed. We illustrate
the four scenarios that arise when removing a rim-hook from
$Y_{p_1,p_2}^{q_1,q_2}$ such that the remaining shape is a hook of
length $2n-c$. We distinguish the case $p_2>0,q_2>0$ (top row),
$p_1=1,p_2=0,q_2>0$ (the case $p_2>0,q_1=1,q_2=0$ being analogous)
(middle row) and $p_1=q_1=1,p_2=q_2=0$ (bottom row).}
\label{F:Young4}
\end{figure}

Next, let us assume $Y\neq Y_{p,q}$. Since $\tau_c$ is a product of
two disjoint cycles, $\chi^Y(\tau_c)\neq 0$ implies that at most two
removals of rim-hooks exhaust $Y$ by the Murnaghan-Nakayama rule, so
$Y$ is necessarily of the form $Y_{p_1,p_2}^{q_1,q_2}$.  Since the
cycles of $\tau_c$ are disjoint, we may remove them in any desired
order, and we choose to first remove a rim-hook of size $c$ and
second a rim-hook of size $2n-c$.  There are exactly four scenarios
for such removals as illustrated in Figure \ref{F:Young4}, it is
worth mentioning that these scenarios can be distinguished by
whether the first rim hook contains no squares in the leftmost
column or the top row (a), contains squares in one but not both of
them (b or c), or contains squares in both (d). We accordingly
derive
\begin{eqnarray*}
\sum_{c=1}^{2n-1}\chi^{Y_{p_1,p_2}^{q_1,q_2}}(\tau_c) & = &
\sum_{c=1}^{2n-1} \left[{(-1)^{p_2}\delta_{c,p_2+q_2+1} (-1)^{p_1}
         \delta_{2n-c,p_1+q_1+1}}_{} +\right. \\
 & & \qquad {(-1)^{p_1-1} \delta_{c,(p_1+1)+q_2}
     (-1)^{p_2} \delta_{2n-c,p_2+1+q_1}}_{} + \\
 & & \qquad {(-1)^{p_2+1} \delta_{c,p_2+(q_1+1)} (-1)^{p_1}
     \delta_{2n-c,p_1+1+q_2}}_{}+\\
 & & \qquad \left. {(-1)^{p_1} \delta_{c,(p_1+1)+q_1} (-1)^{p_2}
     \delta_{2n-c,p_2+1+q_2}}_{}\right] \\
 & = & 0
\end{eqnarray*}
since the first/third and second/fourth terms pairwise cancel.

Thus, for any Young diagram $Y$ other than a $(p,q)$-hook
the term $\sum_{c}\chi^{Y_{p_1,p_2}^{q_1,q_2}}(\tau_c)$ is trivial,
and Lemma~\ref{L:char} follows.
\end{proof}

\begin{lemma}\label{L:Schur}
For any Young diagram $Y$ that contains $2n$ squares, we have
\begin{equation}\label{E:schur}
\frac{1}{(2n)!}\sum_{\pi\in S_{2n}}N^{\sum_i\pi_i}\chi^Y(\pi)=
s_{Y}(1,\dots,1),
\end{equation}
where $s_Y(x_1,\dots,x_N)$ denotes the Schur-polynomial over $N\ge 2n$
indeterminants.
\end{lemma}
\begin{proof}
Rewriting
\begin{equation*}
N^{\sum_i\pi_i}=\prod_{\pi_i}\left(\sum_{h=1}^N 1^{i}\right)^{\pi_i},
\end{equation*}
as a product of power sums $p_i(x_1,\dots,x_N)=\sum_{h=1}^N
x_h^{i}$, we identify eq.~(\ref{E:schur}) via the Frobenius Theorem
as a particular value of the Schur-polynomial $s_{Y}(x_1,\dots,x_N)$
of $Y$ over $N\ge 2n$ indeterminants:
\begin{equation}\label{E:Schur}
\frac{1}{(2n)!}\sum_{\pi\in S_{2n}} \prod_{\pi_i}
\left(\sum_{h=1}^N 1^{i}\right)^{\pi_i}
\, \chi^{Y}(\pi)=
\frac{1}{(2n)!}\sum_{\pi\in S_{2n}}  \chi^{Y}(\pi)\,
\prod_{\pi_i} p_i(1,\dots,1)^{\pi_i}=
s_{Y}(1,\dots,1).
\end{equation}
\end{proof}


\section{The generating function}


In analogy to the polynomials $P(n,x)$ in eq.~(\ref{E:polynN}), we define
\begin{equation}
Q(n,x)=\sum_{\{g\mid  2g\le n\}} \mathbf{c}^{[2]}_g(n) \cdot x^{n-2g}.
\end{equation}

\begin{lemma}\label{L:Q-P}
The polynomial $Q(n,N)$ can be written as
\begin{equation}
Q(n,N)=U(n,N)-V(n,N)
\end{equation}
where
\begin{equation}\label{E:w}
 -3-4Nz-2N^2z^2+2\sum_{n\geq1}\frac{U(n,N)}{(2n-1)!!}z^{n+2}
=  (z+z^3)\frac{d}{dz}\left(\left(\frac{1+z}{1-z}\right)^N\right)-
3\left(\frac{1+z}{1-z}\right)^N
\end{equation}
and
\begin{equation}
V(n,N)=\sum_{1\le c\le n-1}P(c,N)P(n-c,N).
\end{equation}
\end{lemma}
\begin{proof}
A connected chord diagram of genus $g$ on two backbones with $n$ chords
is described by permutations
 $\tau_c=(1,\dots,c)\circ (c+1,\dots,2n)$ and $\iota\in[2^n]$ via Proposition
 \ref{P:isos} satisfies $2-2g-r=2-n$, where the number $r=n-2g$
 of boundary components is given by the
the number of cycles of $\tau_c\circ\iota$.

On the other hand, if the chord diagram corresponding to
$\tau_c$ and $\iota$ is disconnected, then not only
does $\tau_c$ decompose into disjoint cycles $\tau_1=(1,\ldots, c)$
and $\tau_2=(c+1,\ldots ,2n)$ but also
$\iota=\iota_1\circ\iota _2$ similarly
decomposes into disjoint permutations, and
the number of boundary components is given by
\begin{equation}
\sum_i(\tau_c\circ\iota)_i=\sum_s(\tau_1\circ \iota_1)_s+\sum_t(\tau_2\circ\iota_2)_t,
\end{equation}
where $(\pi)_i$ denotes the number of cycles of length $i$ in $\pi$,
$\pi\in S_{2n}$.

We proceed by writing $Q(n,N)=U(n,N)-V(n,N)$ as the difference of two terms,
the first being the contribution of all chord diagrams, irrespective of
being connected, and the second being the contribution of all
disconnected chord diagrams. Thus,
\begin{eqnarray*}
Q(n,N) & = &
        \sum_{\{g \, \mid\,  2g\le n\}} \mathbf{c}^{[2]}_g(n) \, N^{n-2g} \\
         & = &   \underbrace{\sum_{1\le c\le 2n-1}
               \sum_{\iota\in [2^n]\atop \tau_c=(1,\dots,c)(c+1,\dots,2n)}
               N^{\sum_i(\tau_c\iota)_i}}_{U(n,N)}  \\
        & - &  \underbrace{\sum_{\atop 1\le d\le n-1}
               \left(\sum_{\iota_1\in [2^d]\atop \tau_1=(1,\dots,2d)}
               N^{\sum_i(\tau_1\iota_1)_i}\right)
               \left(\sum_{\iota_2\in [2^{n-d}]\atop \tau_2=(2d+1,\dots,2n)}
                N^{\sum_i(\tau_2\iota_2)_i}\right)}_{V(n,N)}
\end{eqnarray*}
since the number of vertices in a chord diagram is necessarily even.
In view of
\begin{equation*}
P(d,N)    =  \sum_{\iota_1\in [2^{d}]\atop \tau_1=(1,\dots,2d)}
                   N^{\sum_i(\tau_1\iota_1)_i}\quad\text{\rm and}\quad
P(n-d,N)  =  \sum_{\iota_2\in [2^{n-d}]\atop \tau_2=(2d+1,\dots,2n)}
                   N^{\sum_i(\tau_2\iota_2)_i},
\end{equation*}
we conclude
\begin{eqnarray*}
V(n,N) & = & \sum_{1\le d\le n-1}P(d,N)P(n-d,N).
\end{eqnarray*}
Turning our attention now to $U(n,N)$, we have
\begin{eqnarray}\label{E:wq}
U(n,N)=
\sum_{c=1}^{2n-1}\sum_{\iota\in [2^n]\atop \tau_c=(1,\dots,c)(c+1,\dots,2n)}
               N^{\sum_i(\tau_c\iota)_i} & = &
              \sum_{[\pi]} N^{\sum_i\pi_i}\sum_{c=1}^{2n-1}
              \sum_{\iota\in [2^n]\atop \tau_c\iota\in [\pi]} 1.
\end{eqnarray}
Expressing the right-hand side of eq.~(\ref{E:wq}) via Kronecker
delta-functions, we obtain a sum taken over all permutations
\begin{eqnarray*}
\sum_{[\pi]} N^{\sum_i\pi_i}\sum_{c=1}^{2n-1}
              \sum_{\iota\in [2^n]\atop \tau_c\iota\in [\pi]} 1
&=& \sum_{[\pi]} N^{\sum_i\pi_i} \sum_{c=1}^{2n-1}\sum_{\sigma\in S_{2n}}
\delta_{[\sigma],[2^n]}\cdot\delta_{[\tau_c\sigma],[\pi]},
\end{eqnarray*}
and application of Lemma~\ref{L:kro} gives
\begin{eqnarray}\label{E:erni}
U(n,N) & = & (2n-1)!!\cdot
\sum_{c=1}^{2n-1}\sum_{Y} \frac{\chi^Y([2^n])\chi^Y(\tau_c)}{\chi^Y([1^{2n}])}
\frac{1}{(2n)!}\sum_{\pi\in S_{2n}}N^{\sum_i\pi_i}
\chi^Y(\pi).
\end{eqnarray}
Interchanging the order of summations and using Lemma~\ref{L:Schur}, we
may rewrite this as
\begin{eqnarray*}\label{E:erni-pq}
U(n,N) & = & (2n-1)!!
\sum_{Y} \frac{\chi^Y([2^n])}{\chi^Y([1^{2n}])}
\left(\sum_{c=1}^{2n-1}\chi^Y(\tau_c)\right)
s_Y(1,\ldots,1).
\end{eqnarray*}
Now, according to Lemma~\ref{L:char}, we have $\sum_{c}\chi^Y(\tau_c)=
(-1)^p(q-p)\,\delta_{Y,Y_{p,q}}$. This reduces the character sum
to the consideration of characters $\chi^{p,q}$ and Schur-polynomials
$s_{p,q}=s_{Y_{p,q}}$ associated to the irreducible representations $Y_{p,q}$.
The corresponding expressions computable
using the Murnaghan-Nakayama rule, for instance,
are given by
\begin{eqnarray}
\chi^{p,q}([1^{2n}]) & = & \binom{2n-1}{q},\\
\chi^{p,q}([2^n]) &= &
\begin{cases}
(-1)^{\frac{p}{2}}\binom{n-1}{\frac{p}{2}}; & \text{\rm for $p\equiv 0
\mod 2$,}\\
(-1)^{\frac{p+1}{2}}\binom{n-1}{\frac{p-1}{2}}; & \text{\rm for $p\equiv 1
\mod 2$,}
\end{cases}\\
s_{p,q}(1,\dots,1)&=&\binom{N+q}{2n}\,\binom{2n-1}{q}.
\end{eqnarray}
Consequently, we compute
\begin{eqnarray*}
\frac{U(n,N)}{(2n-1)!!} & =& \sum_{j=0}^{n-1}(-1)^j\binom{n-1}{j}\left[
(2n-4j-1)\binom{N+2n-2j-1}{2n}\right.\\
&&\qquad\qquad\qquad\qquad\left.+(2n-4j-3)\binom{N+2n-2j-2}{2n}\right] \\
& =& \sum_{j=0}^{n-1}(-1)^j \binom{n-1}{j}\frac{1}{2\pi i}
\oint \left[(2n-4j-1)\frac{(1+x)^{N+2n-2j-1}}{x^{N-2j}}\right.\\
&&\qquad\qquad\qquad\qquad\left.+
(2n-4j-3)\frac{(1+x)^{N+2n-2j-2}}{x^{N-2j-1}}\right]dx.
\end{eqnarray*}
Taking the summation into the integral, we obtain
\begin{eqnarray*}
& =& \frac{1}{2\pi i}
\oint \frac{(1+x)^N}{x^N}
\left(\sum_{j=0}^{n-1}(-1)^j\binom{n-1}{j}(2n-4j-1)
x^{2j}(1+x)^{2n-2j-1}\right.\\
&&\qquad\qquad\qquad\quad+
\left.\sum_{j=0}^{n-1}(-1)^j\binom{n-1}{j}(2n-4j-3)
x^{2j+1}(1+x)^{2n-2j-2}\right)dx \\
 & = & \frac{1}{2\pi i} \oint \frac{(1+x)^N}{x^N}\left( 1+2x\right)^{n-1}
                                         \left((n-1)(1+2x)^2+n\right)dx
\end{eqnarray*}
using $1+2x=(1+x)^2-x^2$.
It remains to substitute $z=1/(1+2x)$ and compute
\begin{eqnarray*}
\frac{2U(n,N)}{(2n-1)!!}
& =& \frac{1}{2\pi i}
 \oint\left(\frac{1+z}{1-z}\right)^N\left(\frac{n-1}{z^{n+3}}
+\frac{n}{z^{n+1}}\right)dz\\
& =& (n-1)[z^{n+2}]\left(\frac{1+z}{1-z}\right)^N
+n[z^n]\left(\frac{1+z}{1-z}\right)^N\\
& =&[z^{n+1}]\frac{d}{dz}\left(\left(\frac{1+z}{1-z}\right)^N\right)
+[z^{n-1}]\frac{d}{dz}\left(\left(\frac{1+z}{1-z}\right)^N\right)
-3[z^{n+2}]\left(\frac{1+z}{1-z}\right)^N
\end{eqnarray*}
so that
\begin{equation*}
-3-4Nz-2N^2z^2+2\sum_{n\geq1}\frac{U(n,N)}{(2n-1)!!}z^{n+2}
 =  (z+z^3)\frac{d}{dz}\left(\left(\frac{1+z}{1-z}\right)^N\right)-
3\left(\frac{1+z}{1-z}\right)^N
\end{equation*}
completing the proof.
\end{proof}

\begin{theorem}\label{T:main}
For any $g\ge 0$, the generating function $\mathbf{C}_g^{[2]}(z)$ is a rational
function with integer coefficients given by
\begin{eqnarray*}
\mathbf{C}_g^{[2]}(z) & = & \frac{P_g^{[2]}(z)}{(1-4z)^{3g+2}},
\end{eqnarray*}
where $P_g^{[2]}(z)$ is an integral polynomials of degree at most
$(3g+1)$, $P_g^{[2]}(1/4)> 0$ and $[z^h]P_g^{[2]}(z)=0$, for $0\leq h\leq 2g$.
Furthermore, we have
\begin{eqnarray}
P_g^{[2]}(z) & = & z^{-1}P_{g+1}(z)-\sum_{g_1=1}^{g}P_{g_1}(z)P_{g+1-g_1}(z), \\
\label{E:22}
\left[{z^{2g+1}}\right]P^{[2]}_g(z) & = & \mathbf{c}^{[2]}_g(2g+1) =
                  \mathbf{c}_{g+1}(2g+2) = \frac{(4g+4)!}{4^{g+1}(2g+3)!}.
\end{eqnarray}
and the coefficients of ${\bf C}_g^{[2]}(z)$ have the asymptotics
\begin{equation}
[z^n]{\bf C}_g^{[2]}(z)\sim \frac{P_g^{[2]}(\frac{1}{4})}{\Gamma(3g+2)}
n^{3g+1} 4^n.
\end{equation}
\end{theorem}
\begin{proof}
Taking the coefficient of $z^{n+2}$ in eq.~(\ref{E:w}), we find
\begin{eqnarray*}
\frac{U(n,N)}{(2n-1)!!} & = &
[z^{n+2}]\biggl ((z+z^3)\sum_{n\ge 0}(n+1)\frac{P(n,N)}{(2n-1)!!}z^n- 3
\sum_{n\ge 0}\frac{P(n,N)}{(2n-1)!!} z^{n+1} \biggr )\\
& = & (n+2)\frac{P(n+1,N)}{(2(n+1)-1)!!} + n \frac{P(n-1,N)}{(2(n-1)-1)!!} -
3\frac{P(n+1,N)}{(2(n+1)-1)!!},
\end{eqnarray*}
for any
$n\ge 1$, whence
\begin{eqnarray*}
  (2n+1)  ~ U(n,N) & = & (n-1)~ P(n+1,N) + n (2n+1)(2n-1)~P(n-1,N).
\end{eqnarray*}
Substituting $U(n,N)=\sum_{2g\le n}\mathbf{u}_g(n) N^{n-2g}$ and
$P(n,N)= \sum_{2g\le n}\mathbf{c}_g(n) N^{n+1-2g}$, we obtain
\begin{eqnarray*}
  (2n+1) \sum_{2g\le n}\mathbf{u}_g(n) N^{n-2g} & = &
   (n-1) \sum_{2g\le n+1}\mathbf{c}_g(n+1) N^{n+2-2g}  \\
& + &  n(2n+1)(2n-1) \sum_{2g\le n-1}\mathbf{c}_g(n-1) N^{n-2g},
\end{eqnarray*}
so
\begin{eqnarray*}
(2n+1) ~ \mathbf{u}_g(n) & = & (n-1) ~\mathbf{c}_{g+1}(n+1) + n
(2n+1)(2n-1)~ \mathbf{c}_g(n-1).
\end{eqnarray*}
Using $(n+2)\, \mathbf{c}_{g+1}(n+1)=2(2n+1)\,\mathbf{c}_{g+1}(n)+
                              (2n+1)n(2n-1)\,\mathbf{c}_{g}(n-1)$
from the recursion eq.~(\ref{E:recursion}), we have
\begin{eqnarray*}
  \mathbf{u}_g(n) & = & \mathbf{c}_{g+1}(n+1)-2\,\mathbf{c}_{g+1}(n),
\end{eqnarray*}
or equivalently, setting $\mathbf{U}_g(z)=\sum_n\mathbf{u}_g(n)z^n$,
it follows that
\begin{eqnarray}\label{E:ww}
z \mathbf{U}_g(z) & = & (1-2z) \mathbf{C}_{g+1}(z).
\end{eqnarray}
We next consider the term $V(n,N)$ as a polynomial in $N$:
\begin{eqnarray}
V(n,N) & = & \sum_{1\le d\le n-1}
\left(\sum_{g_1}\mathbf{c}_{g_1}(d) N^{d+1-2g_1}\right)
\left(\sum_{g_2}\mathbf{c}_{g_2}(n-d)N^{(n-d)+1-2g_2}\right)\nonumber\\
\label{E:contri}
& = & \sum_{g\ge 0}\sum_{g_1=0}^g\sum_{1\le d\le n-1}
      \mathbf{c}_{g_1}(d) \mathbf{c}_{g-g_1}(n-d)N^{n+2-2g},
\end{eqnarray}
where for $i=1,2$, $\mathbf{c}_{g_i}(a)=0$ if $2g_i>a$.
According to eq.~(\ref{E:contri}) for fixed genus $g$,
the contribution of pairs of chord diagrams, each having one backbone, of
genus $g_1$ and $g_2$ to $N^{n-2g}$ is given by
\begin{eqnarray}\label{E:www}
\sum_{g_1=0}^{g+1}
\sum_{1\le d\le n-1}\mathbf{c}_{g_1}(d)\, \mathbf{c}_{g+1-g_1}(n-d)
& = & [z^n]\sum_{g_1=0}^{g+1} \mathbf{C}^*_{g_1}(z)\mathbf{C}^*_{g+1-g_1}(z),
\end{eqnarray}
where
\begin{equation}
\mathbf{C}^*_{g_1}(z)=
\begin{cases}
{\bf C}_0(z)-1;& \text{\rm for $g_1=0$}, \\
\mathbf{C}_{g_1}(z); & \text{\rm otherwise.}
\end{cases}
\end{equation}
This necessary modification stems from the fact that $1\le d\le n-1$ implies that
the coefficient $\mathbf{c}_0(0)=1$ does not appear.
Using eqs.~(\ref{E:ww}) and~(\ref{E:www}), we obtain
\begin{eqnarray*}
\mathbf{C}^{[2]}_g(z) & = & \mathbf{U}_g(z)-
                \sum_{g_1} \mathbf{C}^*_{g_1}(z)\mathbf{C}^*_{g+1-g_1}(z) \\
                & = & \frac{1-2z}{z}\mathbf{C}_{g+1}(z)-
                   \sum_{g_1} \mathbf{C}^*_{g_1}(z)\mathbf{C}^*_{g+1-g_1}(z).
\end{eqnarray*}
Now,  for $g\ge 1$, we have $\mathbf{C}^*_g(z)=\mathbf{C}_g(z)=
\frac{P_{g}(z)\sqrt{1-4z}}{(1-4z)^{3g}}$ so
\begin{eqnarray*}
     \frac{1-2z}{z}\mathbf{C}_{g+1}(z)-2\mathbf{C}^*_{0}(z)\mathbf{C}_{g+1}(z)
      &=& \left[\frac{1-2z}{z}-2\left(\frac{1-\sqrt{1-4z}}{2z}-1\right)\right]
\mathbf{C}_{g+1}(z)\\
&=& \frac{P_{g+1}(z)}{z(1-4z)^{3g+2}},
\end{eqnarray*}
and hence the two non-rational terms conveniently cancel.
We continue by computing
\begin{eqnarray*}
\mathbf{C}_{g_1}(z)\mathbf{C}_{g+1-g_1}(z) =
\frac{P_{g_1}(z)\sqrt{1-4z}}{(1-4z)^{3g_1}}
\frac{P_{g+1-g_1}(z)\sqrt{1-4z}}{(1-4z)^{3(g+1-g_1)}} =
\frac{P_{g_1}(z)P_{g+1-g_1}(z)}{(1-4z)^{3g+2}},~{\rm for}~g_1\ge 1,
\end{eqnarray*}
so all other terms in the sum are also rational expressions.

Thus,
\begin{eqnarray}\label{E:here}
\mathbf{C}_g^{[2]}(z) & = &  \frac{z^{-1}P_{g+1}(z)}{(1-4z)^{3g+2}}-
\sum_{g_1=1}^{g}\frac{P_{g_1}(z)P_{g+1-g_1}(z)}{(1-4z)^{3g+2}},
\end{eqnarray}
and hence
\begin{equation}\label{E:pg2z}
P^{[2]}_g(z)={\bf C}_g^{[2]}(z)(1-4z)^{3g+2}=z^{-1}P_{g+1}(z)-\sum_{g_1=1}^{g}P_{g_1}(z)P_{g+1-g_1}(z)\end{equation}
as was asserted.  Since the $P_{g}(z)$ are polynomials of degree at most $(3g-1)$,
it follows from eq.~(\ref{E:pg2z}) that the degree of $P_g^{[2]}(z)$ is at most
$3g+1$.
Moreover, it follows immediately from $P_g^{[2]}(z)={\bf C}_g^{[2]}(z) (1-4z)^{3g+2}$
that the polynomial $P_g^{[2]}(z)$ has integer coefficients.

We next show that $P^{[2]}_g(1/4)>0$.
According to Lemma~\ref{L:gamma},
$P_{g_1}(1/4)P_{g+1-g_1}(1/4)$ is given by
\begin{eqnarray*}
& &\frac {\Gamma  \left( g_1-\frac{1}{6} \right)
\Gamma \left( g_1+\frac{1}{2}\right)
\Gamma\left( g_1+\frac{1}{6} \right)
\Gamma  \left( g-g_1+\frac{5}{6} \right)
\Gamma  \left( g-g_1+\frac{3}{2}\right)
\Gamma  \left( g-g_1+\frac{7}{6} \right)}
{36\,{\pi }^{3}\,{4}^{g+1}\,{9}^{-(g+1)}\,
\Gamma  \left( g_1+1 \right) \Gamma  \left(g+1- g_1+1 \right) }\\
&\leq&\frac{{9}^{(g+1)}}{36\,{\pi }^{3}\,{4}^{g+1}}
\left(\Gamma \left( g_1-\frac{1}{6} \right)
\Gamma\left( g_1+\frac{1}{6}\right) \Gamma\left( g-g_1+\frac{5}{6}\right)
\Gamma\left( g-g_1+\frac{7}{6} \right)\right),
\end{eqnarray*}
where we use the identity $\Gamma  \left(g+1/2\right)<\Gamma \left(g+1\right)$, which
follows from $\Gamma(x+1)=x\Gamma (x)$. Furthermore
\begin{eqnarray*}
Z_{g_1}&=&\Gamma  \left( g_1-\frac{1}{6} \right)
\Gamma\left( g_1+\frac{1}{6} \right)\Gamma\left(g-g_1+\frac{5}{6}\right)
\Gamma\left( g-g_1+\frac{7}{6} \right)\\
&=&\left(g_1-\frac{7}{6}\right)_{g_1-1}
\left(g_1-\frac{5}{6}\right)_{g_1}
\left(g-g_1-\frac{1}{6}\right)_{g-g_1}
\left(g-g_1+\frac{1}{6}\right)_{g+1-g_1}
\left(\Gamma\left(\frac{5}{6}\right)
\Gamma\left(\frac{1}{6}\right)\right)^2,
\end{eqnarray*}
where
$(x)_n=x(x-1)\cdots(x-n+1)$ denotes the Pochhammer symbol.
It follows that
$Z_{g_1}=Z_{g+1-g_1}$, for $1\leq g_1\leq g$,
and furthermore,
\begin{equation}
Z_{g_1}\leq Z_1=(g-7/6)_{g-1}(g-5/6)_{g}
\left(\Gamma(5/6)\Gamma(1/6)\right)^2/6.
\end{equation}
Thus,
\begin{equation}
\sum_{g_1=1}^gP_{g_1}(1/4)P_{g+1-g_1}(1/4)\leq
\frac{{9}^{(g+1)}}{36\,{\pi }^{3}\,{4}^{g+1}}\cdot\frac{g}{6}\cdot
(g-7/6)_{g-1}(g-5/6)_{g}\left(\Gamma(5/6)\Gamma(1/6)\right)^2.
\end{equation}
We proceed by estimating
\begin{eqnarray*}
4P_{g+1}(1/4) &=&
\frac {{9}^{g+1}}{6\,{\pi }^{3/2}\,{4}^{g}}\cdot\frac{\Gamma
\left( g+5/6 \right) \Gamma  \left( g+3/2\right) \Gamma
\left( g+7/6 \right) }{\Gamma  \left( g+2 \right) }\\
&\geq &\frac {{9}^{g+1}}{6\,{\pi }^{3/2}\,{4}^{g}}\cdot
\frac{\Gamma  \left( g+5/6 \right) \Gamma  \left( g+7/6 \right) }{ g+1  },\\
&= &\frac {{9}^{g+1}}{6\,{\pi }^{3/2}\,{4}^{g}}\cdot
\frac{\left(g-1/6\right)_{g}\left(g+1/6\right)_{g+1}
\Gamma(5/6)\Gamma(1/6)}{g+1}
\end{eqnarray*}
since $\Gamma\left(g+3/2\right)\geq\Gamma \left( g+1\right)$. Thus,
\begin{eqnarray*}
P^{[2]}_g(1/4)&=&4P_{g+1}(1/4)-
\sum_{g_1=1}^gP_{g_1}(1/4)P_{g+1-g_1}(1/4)\\
&\geq& \frac {{9}^{g+1}}{6\,{\pi }^{3/2}\,{4}^{g+1}}
\Gamma(5/6)\Gamma(1/6)(g-7/6)_{g-1}(g-5/6)_{g}\\
&&\quad\times
\left(\frac{4\left(g-1/6\right)
\left(g+1/6\right)}{ g+1  }-
                   \frac{g\,\Gamma(5/6)\Gamma(1/6)}{36\,{\pi }^{3/2}}\right).
\end{eqnarray*}
Finally,
\begin{equation}
\frac{4\left(g-1/6\right)
\left(g+1/6\right)}{ g+1  }-\frac{g\,\Gamma(5/6)\Gamma(1/6)}{36\,{\pi }^{3/2}}
=\frac{4g^2-1/9-g^2/(18\sqrt{\pi})-g/(18\sqrt{\pi})}{g+1}>0,
\end{equation}
for $g\geq 1$, so we indeed have $P^{[2]}_g(1/4)>0$ as claimed.

To see that $[z^h]P_g^{[2]}(z)=0$, for any $0\leq h\leq 2g$, it follows from eq.~(\ref{E:here}) that
\begin{equation}
[z^h]P_g^{[2]}(z)= [z^{h+1}]P_{g+1}(z)-
                   [z^h]\sum_{g_1=1}^gP_{g_1}(z)P_{g+1-g_1}(z).
\end{equation}
Since $[z^h]P_{g_1}=0$, for $h<2g_1$ or $h>3g_1-1$,
we conclude from
$[z^h]\sum_{g_1=1}^gP_{g_1}(z)P_{g+1-g_1}(z)=\sum_{g_1=1}^g
\sum_{i=0}^h[z^i]P_{g_1}(z)[z^{h-i}]P_{g+1-g_1}(z)$ that
$[z^i]P_{g_1}(z)[z^{h-i}]P_{g+1-g_1}(z)\neq 0$ implies
$i\geq 2g_1$ and $h-i\geq 2(g+1-g_1)$.
Thus,
\begin{equation}\label{E:-1}
[z^h]\sum_{g_1=1}^gP_{g_1}(z)P_{g+1-g_1}(z)=0,
\end{equation}
for $0\leq h\leq 2g+1$, and consequently,
\begin{equation}\label{E:0}
[z^h]P_g^{[2]}(z)=[z^{h+1}]P_{g+1}(z)=0,
\end{equation}
for $0\leq h\leq 2g$, as required.

It remains only to compute the coefficient of $z^{2g+1}$ in $P_g^{[2]}(z)$.
To this end, we have
\begin{equation}
[z^{2g+1}]P_g^{[2]}(z)=[z^{2g+2}]P_{g+1}(z)-[z^{2g+1}]
\sum_{g_1=1}^gP_{g_1}(z)P_{g+1-g_1}(z)=[z^{2g+2}]P_{g+1}(z)
\end{equation}
in light of eq.~(\ref{E:-1}). Since $[z^h]P_g^{[2]}(z)=0$, for any
$0\leq h\leq 2g$, we conclude from $P_g^{[2]}(z)={\bf C}_g^{[2]}(z)
(1-4z)^{3g+2}$ that
$[z^{2g+1}]P_g^{[2]}(z)=\mathbf{c}^{[2]}_g(2g+1)$, and hence using
eq.~(\ref{E:values})
\begin{equation}
\mathbf{c}^{[2]}_g(2g+1) = [z^{2g+1}]P_g^{[2]}(z)=
[z^{2g+2}]P_{g+1}(z)=\mathbf{c}_{g+1}(2g+2)=\frac{(4g+4)!}{4^{g+1}(2g+3)!}.
\end{equation}
\end{proof}

\begin{corollary}\label{C:formel}
We have the explicit expressions
\begin{eqnarray*}
\mathbf{c}^{[2]}_0(n) & =&  n\, 4^{n-1},\\
\mathbf{c}^{[2]}_1(n) & =&  \frac{1}{12}\left( 13\,n+3 \right)
                             n\left( n-1 \right)\left( n-2 \right)\,{4}^{n-3},\\
\mathbf{c}^{[2]}_2(n) & =&  {\frac {1}{180}}\,\left(
                             445\,{n}^{2}-401\,n-210 \right) n
                       \left(n-1 \right)\left(n-2 \right)\left( n-3 \right)
                                              \left(n-4 \right)\, {4}^{n-6} .
\end{eqnarray*}
\end{corollary}
\begin{proof}
Using the expression for $P_1(z)$ in the Introduction, Theorem~\ref{T:main}
gives
\begin{equation}
\mathbf{C}^{[2]}_0(z)=\frac{z^{-1}P_1(z)}{(1-4z)^2}=\frac{z}{(1-4z)^2},
\end{equation}
which immediately implies
\begin{equation}
\mathbf{c}^{[2]}_0(n)=[z^n]\mathbf{C}^{[2]}_0(z)
=[z^{n-1}]\left(\frac{1}{1-4z}\right)^2
=\sum_{i=0}^{n-1}4^{i}4^{n-1-i}=n4^{n-1}.
\end{equation}

In order to derive the expression for $\mathbf{c}^{[2]}_1(n)$, we
use both of the expressions $P_1(z)=z^2$ and $P_2(z)=21z^5+21z^4$.
Theorem~\ref{T:main}
implies
\begin{equation}
\mathbf{C}^{[2]}_1(z)=\frac{z^{-1}P_2(z)-P_1(z)^2}{(1-4z)^5}
=\frac{(20z+21)z^3}{(1-4z)^5},
\end{equation}
and differentiation gives
\begin{equation}\label{E:ode1}
(52z^3-13z^2)\frac{{d}^2\mathbf{C}_1(z)}{{d}z^2}+
(168z^2+36z)\frac{{d}\mathbf{C}_1(z)}{{d}z}+(64z-30)\mathbf{C}_1(z)=0,
\end{equation}
where $\frac{{d}^3\mathbf{C}_1(z)}{{d} z^3}|_{z=0}=126$.
Straightforward calculation shows that eq.~(\ref{E:ode1}) implies
\begin{equation}
\mathbf{c}^{[2]}_1(n+1)={52n^2+116n+64\over 13n^2-23n-6}
\mathbf{c}^{[2]}_1(n),
\end{equation}
where $\mathbf{c}^{[2]}_1(i)=0$, for $0\leq i\leq 2$ and
$\mathbf{c}^{[2]}_1(3)=21$. It remains to observe that
$\mathbf{c}^{[2]}_1(n)=\frac{1}{12} \left( 13\,n+3 \right)  n\left(
n-1 \right)\left( n-2 \right)\,{4}^{n-3}$ satisfies this recursion
together with its initial condition.

To compute $\mathbf{c}^{[2]}_2(n)$, we likewise employ
$P_3(z)=11z^6
\left( 158\,{z}^{2}+558\,z+135 \right)$ and Theorem~\ref{T:main} to conclude
\begin{equation}
\mathbf{C}^{[2]}_2(z)=\frac{z^{-1}P_3(z)-2P_1(z)P_2(z)}{(1-4z)^8}
=\frac{(1696z^2+6096z+1485)z^5}{(1-4z)^8}.
\end{equation}
This yields the ODE
\begin{eqnarray*}
&&(1780z^4-445z^3)\frac{{d}^3\mathbf{C}^{[2]}_2(z)}{{d}z^3}+
(9076z^3+2181z^2)\frac{{d}^2\mathbf{C}^{[2]}_2(z)}{{d} z^2}\\
&&+(6808z^2-4020z)\frac{{d}\mathbf{C}^{[2]}_2(z)}{{d} z}-(664z-3180)
\mathbf{C}^{[2]}_2(z)=0,
\end{eqnarray*}
where $\frac{{d}^5\mathbf{C}^{[2]}_2(z)}{{d} z^5}|_{z=0}=178200$.
Thus,
\begin{equation}
\mathbf{c}^{[2]}_2(n+1)={(1780n^3+3736n^2+1292n-664)
\over(445n^3-2181n^2+1394n+840)}\mathbf{c}^{[2]}_2(n),
\end{equation}
where $\mathbf{c}^{[2]}_2(i)=0$, for $0\leq i\leq 4$, with
$\mathbf{c}^{[2]}_2(5)=1485$, and we conclude as before
the verity of the asserted formula.
\end{proof}

\end{document}